\documentclass[12pt,a4paper]{article}
\usepackage{amssymb,enumerate,amsmath}
\usepackage{graphicx}
\usepackage[latin1]{inputenc}
\usepackage{pgf,tikz}\usetikzlibrary{arrows}
\usepackage{esvect}\renewcommand{\vec}[1]{\vv{#1}}	
\newtheorem{lm}{{\bf Lemma}}
\newtheorem{theor}[lm]{{\bf Theorem}}
\newtheorem{deff}[lm]{{\bf Definition}}
\newtheorem{exemple}[lm]{{\bf Example}}
\newcommand{\sep}{\hspace{-2mm}{\bf.}}
\newcommand{\rf}[1]{(\ref{#1})}
\newcommand{\df}[2]{\begin{deff}{\label{#1}}\sep{\rm #2}\end{deff}}
\newcommand{\ex}[2]{\begin{exemple}{\label{#1}}\sep{\rm
#2}\end{exemple}}
\newcommand{\lem}[2]{\begin{lm}\label{#1}\sep\it#2\end{lm}}
\newcommand{\theo}[2]{\begin{theor}{\label{#1}}\sep{\it
#2}\end{theor}}
\newcommand{\eq}[2]{\begin{equation}#2\ \ \  \ \label{#1}\end{equation}}
\newcommand{\pr}{\noindent {\it Proof. }}
\newcommand{\ep}{\hfill$\square$\medskip}
\newcommand{\tq}{\ ;\ }
\newcommand{\be}{\begin{enumerate}[(i)]}
\newcommand{\ee}{\end{enumerate}}
\def\bd{{{\rm bd\,}}}
\def\inn{{\rm int\,}}
\renewcommand{\bar}{\overline}
\newcommand{\tl}{\to_{\ell}\,}
\newcommand{\tr}{\to_{r}\,}
\renewcommand{\d}{\delta}
\newcommand{\eps}{\varepsilon}
\newcommand{\g}{\gamma}
\newcommand{\q}{\theta}
\renewcommand{\S}{{\mathbb S}}
\newcommand{\V}{{\mathcal V}}
\newcommand{\N}{\mathbb{N}}

\newcommand{\R}{\mathbb{R}}
\newcommand{\id}{{\bf id}}
\newcommand{\is}{{\rm Isom^+\R^2}\,}
\title{
{\Large Fixing and almost fixing a planar convex body}}
\author{Augustin\ Fruchard}
\date{\today}
\begin{document}
\maketitle
\begin{abstract}
A set of points $a_1,\dots,a_n$ {\sl fixes} a planar convex body $K$
if the points are on $\bd K$, the boundary of $K$, and if 
any small move of $K$ brings some
point of the set in $\inn K$, the interior of $K$.
The points $a_1,\dots,a_n\in\bd K$ {\it almost fix} $K$ if,
for any neighbourhoods $\V_i$ of $a_i$ $(i=1,\dots,n)$,
there are pairs of points $a_i',a_i''\in\V_i\cap\bd K$
such that $a_1',a_1'',\dots,a_n''$ fix $K$.
This note compares several definitions of these notions and
gives first order conditions for $a_1,\dots,a_n\in\bd K$ to fix,
and to almost fix, $K$.
\end{abstract}
\medskip

\medskip

\

\noindent{\bf Notation.} \
The plane $\R^2$ is endowed with its euclidean norm 
$\|\;\|$. The origin of $\R^2$ is denoted by $\bf0$.
The unit cercle of centre $\bf0$ and radius $1$ is denoted by $\S_1$.
The group of planar affine rotations and translations is denoted by $\is$.
The identity of $\R^2$ is denoted by $\id$.
The rotation of centre $\omega\in\R^2$ and angle $\alpha\in\R$ is denoted by
$R_{\omega,\alpha}$.
Given an oriented line $D$, the open half-plane on the left of $D$ is
denoted by $D^+$.
As usual, for $M\subset\R^2$, $\bd M$ denotes its
boundary and $\inn M$ denotes its interior.

Throughout this note, $K$ is a convex compact subset of $\R^2$ of
nonempty interior, and $a_1,\dots,a_n$ are points on $\bd K$.
\df{d1}{
The points $a_1,\dots,a_n$ {\sl fix} $K$ if there is a 
neighbourhood $\V$ of  $\id$ in $\is$ such that, 
for every $f\in\V$ satisfying $f(K)\ne K$,
at least one of the $a_i$ belongs to $\inn f(K)$.
}
This is slightly different from the commonly used definition of fixing points,
e.g. in~\cite{BMU,BFMM}, where a subset $H$ of $\bd K$ is said to fix $K$ if
$\id$ is isolated in the set
$\big\{f\in\is\tq H\cap f(\inn K)=\emptyset\big\}$.
In particular, we consider that a rotation of a disc around its centre
does not move the disc.
If $K$ is not a disc, then both definitions are equivalent,
since in this case there is a neighbourhood $\V$ of $\id$ in $\is$
such that every $f\in\V\setminus\{\id\}$ satisfies $f(K)\ne K$.

Although  most convex bodies in the sense of Baire categories can be fixed by
three points, there are convex bodies which cannot be fixed by any
set of three points, e.g. parallelograms. 
Other examples can be found in~\cite{CSU}.

Another commonly used definition of fixing points, closer to the intuition,
is the following one, see e.g.~\cite{CSU}.
\df{d2}{
The points $a_1,\dots,a_n\in\bd K$ {\sl weakly fix $K$} if,
for any path $\g:[0,1]\to\is,\;t\mapsto\g_t$ such that
$\g_0(K)=K$ and $\g_1(K)\ne K$, there exist $i\in\{1,\dots,n\}$ and
$t\in[0,1]$ such that $a_i\in\inn\g_t(K)$.
}
Obviously, if $a_1,\dots,a_n$ fix $K$ in the sense of Definition~\ref{d1},
then they weakly fix $K$.
The converse is false, as shows the following example.
\ex{e1}{
Choose $K$ admitting $\S_1$ as an inscribed circle and such that
$\bd K\cap\S_1$ is the union of an arc of circle of
length greater than $\pi$ and a countable number of points accumulating
to some point of $\S_1$, say:
$$
\bd K\cap\S_1=\{(\cos\q,\sin\q)\tq
|\q|\in[\pi/4,\pi]\mbox{ or }\exists n\in\N\tq|\q|=\pi/2^n\}.
$$
This can be done by completing the arc of circle of length $3\pi/2$
by  segments tangent to $\S_1$ at these points 
$(\cos(\pi/2^n),\sin(\pm\pi/2^n))$; this could also be done such that $K$
is strictly convex and smooth.
Then the four points $a_i=(\cos(i\pi/2),\sin(i\pi/2))$, $i=1,2,3,4$,
weakly fix $K$ but do not fix it in the sense of Definition~\ref{d1}.
}
{\noindent\bf A  first order condition to fix $K$.} \
\medskip

We want to give an almost necessary and sufficient condition
of first order for $a_1,\dots,a_n$ to fix $K$.
``Of first order'' means that these conditions make use of
the left and right tangents of $\bd K$ at the $a_i$s.
``Almost'' means that the necessary condition involves non-strict
inequalities, whereas the sufficient condition involves the
corresponding strict inequalities.
In the case of three points, we find again the well-known
necessary condition that the three normals are concurrent.
In~\cite{BMU} an almost necessary and sufficient condition 
{\it of second order}, i.e. using the curvature of $\bd K$, is given
such that three points fix a ${\mathcal C}^2$ convex body $K$.
\medskip

We first introduce some notation. We orient $\bd K$ counterclockwise.
Given $a\in\bd K$, let $T_{\ell}(a)$, resp. $T_r(a)$, be the left, resp. 
right, tangent at $\bd K$ in $a$. We orient these lines as $\bd K$;
thus we have $K\subset T_{\ell}(a)^+\cap T_r(a)^+$.
Let $N_\ell(a)$, resp. $N_r(a)$, be the left and right normals at 
$\bd K$ in $a$, oriented in the directions $T_{\ell}(a)^+$,
resp. $T_r(a)^+$, i.e. the line orthogonal to $T_\ell(a)$, resp. $T_r(a)$,
containing $a$, and pointing (as) inward $K$ (as possible).
Let $L(a)$ be the open sector, union of the left open
half-planes bounded by $N_\ell(a)$ and $N_r(a)$:
$$
L(a)=N_\ell(a)^+\cup N_r(a)^+.
$$
Let $\bar{L}(a)$ be the corresponding closed sector,
and let $\vec{L}(a)=\{x\in\mathbb{S}_1\tq a+x\in\bar{L}(a)\}$
be the {\it set of directions} of $\bar{L}(a)$.
Let $R(a)$ and $\bar{R}(a)$ be the analogous 
sectors for the right half-planes.
The set of directions of $\bar{R}(a)$ is $-\vec{L}(a)$, 
hence will not be needed.
If $\bd K$ is differentiable at $a$, then
$\bar{L}(a)=\R^2\setminus{R}(a)$, otherwise $L(a)\cap R(a)$ is the
union of two sectors of vertex $a$.
\theo{t1}{
Let $K$ be a planar convex body and $a_1\dots,a_n\in\bd K$.
\be
\item
If $a_1,...,a_n$ weakly fix $K$, then both intersections
$L(a_1)\cap...\cap L(a_n)$ and $R(a_1)\cap\dots\cap R(a_n)$ are empty.
\item
If the three intersections $\bar{L}(a_1)\cap\dots\cap\bar{L}(a_n)$,
$\bar{R}(a_1)\cap\dots\cap\bar{R}(a_n)$, and
$\vec L(a_1)\cap\dots\cap\vec L(a_n)$ are empty,
then $a_1,\dots,a_n$ fix $K$.
\ee
}
{\noindent\bf Remark.} 
The set of directions is needed in item (ii).
For instance, if $K=[-2,2]\times[0,1]$,
then the three points $a=(-1,0)$, $b=(1,0)$, $c=(0,1)$ do not
fix $K$, although the intersections $\bar L(a)\cap\bar L(b)\cap\bar L(c)$
and $\bar R(a)\cap\bar R(b)\cap\bar R(c)$ are empty.
\medskip

For the proof of Theorem~\ref{t1}, we will need the following statement.
Recall that  $R_{\omega,\alpha}$ is the rotation of centre $\omega$ and
angle $\alpha$.
\lem{l3.1}{
Let $K$ be a planar convex body, let $a\in\bd K$, and let $\omega\in\R^2$. 
\be
\item
If $\omega\in L(a)$ then there exists $\eps>0$ such
that $R_{\omega,-\alpha}(a)\notin K$ for all $\alpha\in\,]0,\eps[$.
\item
Let $(\alpha_\nu)_{\nu\in\N}$ be a sequence of positive numbers tending to $0$
and $(\omega_\nu)_{\nu\in\N}$ be a sequence of points of $\R^2$ tending to
$\omega$.
If $R_{\omega_\nu,-\alpha_\nu}(a)\notin\inn K$ for all $\nu\in\N$,
then $\omega\in\bar{L}(a)$.
\item
Let $(\alpha_\nu)_{\nu\in\N}$ be a sequence of positive numbers tending to $0$
and $(\omega_\nu)_{\nu\in\N}$ be a sequence of points of $\R^2$ tending to
infinity (i.e. $\|\omega_\nu\|\to+\infty$) such that the sequence 
$\Big(\frac{\omega_\nu}{\|\omega_\nu\|}\Big)_{\nu\in\N}$ converges to some
$x\in\S^1$. We assume that the sequence of rotations
$(R_{\omega_\nu,-\alpha_\nu})_{\nu\in\N}$ tends to $\id\in\is$, i.e.
$\displaystyle\lim_{\nu\to+\infty}\alpha_\nu\|\omega_\nu\|=0$.
If $R_{\omega_\nu,-\alpha_\nu}(a)\notin\inn K$  for all $\nu\in\N$,
then $x\in\vec{L}(a)$.
\ee
}
\pr {\it(i)} \ 
If $\omega\in L(a)\cap T_\ell(a)^+$ then one can choose $\eps=2\angle na\omega$,
where $n$ is any point on $N_\ell(a)\cap T_\ell(a)^+$.
If $\omega\in L(a)\setminus T_\ell(a)^+$, then $\eps=\pi$ suits.
\medskip

{\noindent\it(ii)} \
By contraposition, let $\omega\notin\bar{L}(a)$ be fixed, at a distance $2d>0$
from $\bar{L}(a)$ and at a distance $D/2$ from $a$.
Then, for $\nu$ large enough, $\omega_\nu$ is at a distance at least $d$ from 
$\bar{L}(a)$ and at most $D$ from $a$, hence $a\omega_\nu$ makes an angle at
least $\d=\arctan\tfrac dD$ with $N_\ell(a)\cap T_\ell(a)^+$ and with
$N_r(a)\setminus T_r(a)^+$. Let $T'_\ell(\d)=R_{a,-\d}(T_\ell(a))$ and
$T'_r(\d)=R_{a,\d}(T_r(a))$.
By definition of the left and right tangents, there exists $r>0$
such that $D(a,r)\cap\big(T'_\ell(\d)\big)^+\cap\big(T'_r(\d)\big)^+$ is
entirely included in $K$, where $D(a,r)$ denotes the disc of centre $a$
and radius $r$. Then, as soon as $\alpha_\nu<\tfrac rD$, we have
$R_{\omega_\nu,-\alpha_\nu}(a)\in\inn K$.
\medskip

{\noindent\it(iii)} \ 
By contraposition, assume $x\notin\vec L(a)$. Let $\d>0$ be so small
that $x\notin R_{{\bf0},-2\d}\big(\vec L(a)\big)
\cup R_{{\bf0},2\d}\big(\vec L(a)\big)$.
Then, for $\nu$ large enough, we have $\frac{\omega_\nu}{\|\omega_\nu\|}\notin
R_{{\bf0},-\d}\big(\vec L(a)\big)\cup R_{{\bf0},\d}\big(\vec L(a)\big)$.
With the above notation, let $r>0$ be such that
$D(a,r)\cap\big(T'_\ell(\d)\big)^+\cap\big(T'_r(\d)\big)^+\subset K$.
Then, as soon as $\alpha_\nu\|\omega_\nu\|<r$, we have
$R_{\omega_\nu,-\alpha_\nu}(a)\in\inn K$.
\ep

{\noindent\sl Proof of Theorem~\ref{t1}.}
{\it(i)} \ 
By contraposition, assume $L(a_1)\cap\dots\cap L(a_n)$ nonempty and choose
a point $\omega$ in this intersection. 
By Lemma~\ref{l3.1} (i), for each $i\in\{1,\dots,n\}$,
since $\omega\in L(a_i)$, there exists $\eps_i>0$ such that
$R_{\omega,-\alpha}(a_i)\notin K$  for all $\alpha\in]0,\eps_i[$.
As a consequence, for all $\alpha<\min(\eps_1,\dots,\eps_n)$, none of the $a_i$
belongs to $\inn R_{\omega,\alpha}(K)$, hence $a_1,\dots,a_n$ do not weakly
fix $K$.
The case $R(a_1)\cap\dots\cap R(a_n)\ne\emptyset$ is similar, using
clockwise rotations.
\medskip

{\it(ii)} \ 
By contraposition, if $a_1,\dots,a_n$ do not fix $K$ then, for any
neighborhood $\V$ of $\id$ in $\is$, there exists
$f\in\V\setminus\{\id\}$ such that
\eq9{
\forall i=1,\dots,n\quad a_i\notin\inn f(K).
}
It follows that there is a sequence $(f_\nu)_{\nu\in\N}$ 
tending to $\id$ in $\is\setminus\{\id\}$ and satisfying~\rf9.
At least one of the following cases must occur:

- an infinite number of $f_\nu$ are direct (i.e. counterclockwise) rotations,

- an infinite number of $f_\nu$ are indirect rotations,

- an infinite number of $f_\nu$ are translations.

In the first case, considering a subsequence if necessary, we assume without
loss of generality that the whole sequence $(f_\nu)$ consists of direct
rotations, of centres $\omega_\nu\in\R^2$ and angles $\alpha_\nu>0$ tending to $0$. 

- If the sequence $(\omega_\nu)_{\nu\in\N}$ is 
bounded, then a subsequence converges to some point $\omega\in\R^2$. 
Then, Lemma~\ref{l3.1} (ii) and~\rf9 imply that $\omega$ belongs to 
$\bar{L}_1\cap\dots\cap\bar{L}_n$.

- If the sequence $(\omega_\nu)_{\nu\in\N}$ is 
unbounded, then a subsequence $(\omega_{\nu_k})_{k\in\N}$ tends to infinity.
By compactness of $\S^1$, a subsequence of 
$\Big(\frac{\omega_{\nu_k}}{\|\omega_{\nu_k}\|}\Big)_{k\in\N}$
converges to some $x\in\S^1$.
Then Lemma~\ref{l3.1} (iii) and~\rf9 imply that 
$x\in\vec L_1\cap\dots\cap\vec L_n$.
\medskip

The second case of indirect rotations yields similarly
$\bar R_1\cap\dots\cap\bar R_n\ne\emptyset$
or $\vec L_1\cap\dots\cap\vec L_n\ne\emptyset$.
\medskip

In the third case, we also assume without loss of generality that every
$f_\nu$ is a translation of vector $v_\nu$. Then
a subsequence of $\Big(\frac{v_\nu}{\|v_\nu\|}\Big)_{\nu\in\N}$ converges
to some $x\in\S^1$ and we obtain similarly to the proof of Lemma~\ref{l3.1}
(iii) that $R_{{\bf0},-\pi/2}(x)\in\vec L_1\cap\dots\cap\vec L_n$.
\ep
\bigskip

{\noindent\bf A first order condition to almost fix $K$.}
\df{d3}{
The points $a_1,\dots,a_n$ {\it almost fix} $K$ if, 
for any neighbourhoods $\V_i$ of $a_i$ $(i=1,\dots,n)$,
there are pairs of points $a_i',a_i''\in\V_i\cap\bd K$,
such that $a_1',a_1'',\dots,a_n''$ fix $K$.
}
In~\cite{fz} it is proven that any planar convex body can be almost 
fixed by two or three points.
\medskip

It seems natural that, if $a_1,\dots,a_n$ almost fix $K$,
then one can choose $a_i'$ and $a_i''$ on each side of $a_i$.
The following example shows that this is not the case!
\ex{e2}{
For any integer $n\ge2$, let
$p_n=\tfrac1{\cos(\pi/4^n)}\big(\cos(3\pi/4^n),\sin(3\pi/4^n)\big)$
and choose for $K$ the convex hull of $\S_1\cup\{p_2,p_3,\dots\}$.
Then $K$ admits $\S_1$ as an inscribed circle and, for any $\q\in\R$,
we have 
$$(\cos\q,\sin\q)\in\bd K\cap\S_1
\Leftrightarrow
\big(\q\in[\pi/4,2\pi]\mbox{ or }\exists n\in\N\tq
\pi/4^n\le\q\le2\pi/4^n\big).
$$
The rest of $\bd K$ is made of segments tangent to $\S_1$,
crossing at the $p_n$.

Then the four points $a_k=(\cos(k\pi/2),\sin(k\pi/2))$, $k=1,2,3,4$,
almost fix $K$ but in the neighbourhood of $a_4$ one has to choose 
both points $a'_4$ and $a''_4$ on the same side of $a_4$.
}
Our next statement is the analog of Theorem~\ref{t1}.
The only difference is that the large
sectors $L$ and $R$ are replaced by the following small sectors
$\ell$ and $r$: 
$$
\ell(a)=N_\ell(a)^+\cap N_r(a)^+=\R^2\setminus\bar{R}(a),
$$
$r(a)=\R^2\setminus\bar{L}(a)$,
$\bar{\ell}(a)$ and $\bar{r}(a)$ are the corresponding closed sectors
and
$$
\vec{\ell}(a)=\{x\in\mathbb{S}_1\tq a+x\in\bar{\ell}(a)\}.
$$
\theo{t2}{
Let $K$ be a planar convex body and $a_1\dots,a_n\in\bd K$.
\be
\item
If $a_1,...,a_n$ almost fix $K$, then both intersections
$\ell(a_1)\cap...\cap\ell(a_n)$ and $r(a_1)\cap\dots\cap r(a_n)$ are empty.
\item
If the three intersections $\bar{\ell}(a_1)\cap\dots\cap\bar{\ell}(a_n)$,
$\bar{r}(a_1)\cap\dots\cap\bar{r}(a_n)$, and
$\vec\ell(a_1)\cap\dots\cap\vec\ell(a_n)$ are empty,
then $a_1,\dots,a_n$ almost fix K.
\ee
}
\pr
The notation $a'\tl a$ (resp. $a'\tr a$) means that $a'$ tends to $a$ 
from the left (resp. from the right) on $\bd K$, with $a'\ne a$.
The semicontinuity of the left and right tangents can be expressed as
follows:
\eq3{
\lim_{a'\tl a}N_\ell(a')=\lim_{a'\tl a}N_r(a')=N_\ell(a),\quad
\lim_{a''\tr a}N_\ell(a'')=\lim_{a''\tr a}N_r(a'')=N_r(a).
}
Here we use the topology of uniform convergence on every compact set:
A subset $M'$ of $\R^2$ {\it tends to} $M\subset\R^2$ if, for all $r>0$,
the Pompeiu-Hausdorff distance between $M'\cap D({\bf0},r)$ and 
$M\cap D({\bf0},r)$ tends to $0$, where $D({\bf0},r)$ is the closed disc
of centre ${\bf0}$ and radius $r$.
As a consequence of~\rf3, we have 
\begin{align}
\label{8}
\lim_{a'\tl a,a''\tr a}\bar{L}(a')\cap\bar{L}(a'')&=\bar{\ell}(a),\nonumber\\
\lim_{a'\tl a,a''\tr a}\bar{R}(a')\cap\bar{R}(a'')&=\bar{r}(a),\\
\lim_{a'\tl a,a''\tr a}\vec{L}(a')\cap\vec{L}(a'')&=\vec{\ell}(a)\nonumber.
\end{align}

{\it(i)} \ 
Assume $\ell(a_1)\cap...\cap\ell(a_n)$ nonempty and choose a point $\omega$
in this intersection. Let us fix the index $i\in\{1,\dots,n\}$.
By~\rf3, there exists $\eps>0$ 
such that, for all $a'_i,a''_i\in\bd K$
satisfying $|a_ia'_i|<\eps$ and $|a_ia''_i|<\eps$,
we have $\omega\in L(a'_i)\cap L(a''_i)$. By Theorem~\ref{t1} (i), it
follows that $a'_1,\dots,a''_n$ do not fix $K$.
\medskip

{\it(ii)} \ 
By contraposition, if $a_1,\dots,a_n$ do not almost fix $K$,
then there exists $\eps>0$ such that, for all $a'_i,a''_i\in\bd K$
satisfying $|a_ia'_i|<\eps$, one on each side of $a_i$,
$a'_1,\dots,a''_n$ do not fix $K$.
Given sequences $(a'_{i,\nu})_{\nu\in\N}$ and $(a''_{i,\nu})_{\nu\in\N}$
tending to $a_i$, using Theorem~\ref{t1} (ii) and choosing 
adequate subsequences if necessary, it follows that one of the intersections
$\bar{L}(a'_{1,\nu})\cap\dots\cap\bar{L}(a''_{n,\nu})$,  
$\bar{R}(a'_{1,\nu})\cap\dots\cap\bar{R}(a''_{n,\nu})$, or  
$\vec{L}(a'_{1,\nu})\cap\dots\cap\vec{L}(a''_{n,\nu})$, is nonempty
for all $\nu\in\N$.
By~\rf8, we obtain $\bar\ell(a_1)\cap\dots\cap\bar\ell(a_n)$ nonempty in
the first case, 
$\bar{r}(a_1)\cap\dots\cap\bar{r}(a_n)$ nonempty in the second case, and
$\vec\ell(a_1)\cap\dots\cap\vec\ell(a_n)$ nonempty in the last case.
\ep

\end{document}